# Wind Power Providing Flexible Ramp Product

Runze Chen, *Student Member, IEEE*, Jianhui Wang, *Senior Member, IEEE*, Audun Botterud, *Member, IEEE*, Hongbin Sun, *Senior Member,* IEEE

*Abstract*—The deepening penetration of renewables in power systems has contributed to the increasing needs for generation scheduling flexibility. Specifically, for short-term operations, flexibility here indicates that sufficient ramp capacities should be reserved to respond to the expected changes in the load and intermittent generation, also covering a certain amount of their uncertainty. To address the growing requirements for flexible ramp capacity, markets for ramp products have been launched in practice such as the ones in California ISO and Midcontinent ISO. Sometimes, to guarantee sufficient ramp capacity, expensive fast-start units have to be committed in real-time. Moreover, with higher penetration of renewable generation, the flexibility provided by the conventional units might not be enough. Actually, wind power producers are physically capable of offering flexibility, which is sometimes also economically efficient to the entire system. In this paper, we aim to explore the mechanism and possibility of including wind power producers as ramp providers to increase the supply of flexibility. To conduct the analyses, a two-stage stochastic real-time unit commitment model considering ramp capacity adequacy is formulated. Case studies indicate that both the system and the wind power producers can benefit if the wind power is allowed to provide flexible ramp products.

## Nomenclature

### A. Indices, parameters, sets and functions

| | |
|---|---|
| $\omega$ | Penalty factor for flexible ramp capacity (FRC) shortages. |
| $\pi^u, \pi^d$ | Extra compensation factor for upward and downward FRC provided by WPPs. |
| $\xi^u, \xi^d$ | Required upward and downward FRC margin for covering uncertainty. |
| $\rho$ | Probability of each scenario. |
| $\eta^{l,u}, \eta^{l,d}$ | Factors representing the contribution of load forecast to FRC requirement margins. |
| $\eta^{w,u}, \eta^{w,d}$ | Factors representing the contribution of wind forecast to FRC requirement margins. |
| $G^{II}$ | Set of type II conventional units (startup and shutdown processes last for 15 min). |
| $i, G$ | Index and set of units. |
| $j, W$ | Index and set of WPPs. |
| $L$ | System load. |
| $\hat{L}$ | System load forecast. |
| $M$ | A number sufficiently large. |
| $P^{g,max}$ | Maximum output of the units. |
| $P^{g,min}$ | Minimum output of the units. |
| $\tilde{p}^w$ | Offered power output of a WPP. |
| $\hat{p}^w$ | Forecasted power output of a WPP. |
| $\underline{P}^{w,\alpha}$ | Lower quantile of wind power distribution at the confidence level of $\alpha$. |
| $Q^u, Q^d$ | Minimum system upward and downward regulation reserve requirements. |
| $ramp$ | Ramping rate limit of a unit. |
| $ramp^w$ | Ramping rate limit of a WPP. |
| $s, S$ | Index and set of scenarios. |
| $SDRR$ | Maximum shutdown ramping rate of a unit. |
| $SURR$ | Maximum startup ramping rate of a unit. |
| $t, T$ | Index and set of time periods. |
| $T^{15}$ | Set of time periods for unit commitment. |
| $VOLL$ | Value of lost loads. |

### B. Variables

| | |
|---|---|
| $\Delta l$ | Amount of load shedding. |
| $p^g$ | Scheduled power output of a unit. |
| $\bar{p}^{g,su}, \underline{p}^{g,su}$ | Auxiliary variables representing the upper and lower limits of power output of a unit during startup. |
| $\bar{p}^{g,sd}, \underline{p}^{g,sd}$ | Auxiliary variables representing the upper and lower limits of power output of a unit during shutdown. |
| $p^w$ | Scheduled power output of a WPP. |
| $\Delta p^{g,su}, \Delta p^{g,sd}$ | Passive ramping of the units during startup and shutdown. |
| $q^{g,u}, q^{g,d}$ | Upward and downward regulation reserve provided by a unit. |
| $r^{w,u}, r^{w,d}$ | Upward and downward FRC provided by a WPP. |
| $r^{g,u}, r^{g,d}$ | Upward and downward FRC provided by a unit. |
| $\Delta r^u, \Delta r^d$ | Shortage of upward and downward FRCs. |
| $u, v$ | {0,1}, startup and shutdown state variable of a unit. |
| $x$ | {0,1}, on/off state variable of a unit. |

## I. Introduction

OPERATING the system to maintain power balance calls for adequate controllable resources to alter their output to match the varying load. Maintaining this balance can be fairly challenging when facing the deepening penetration of variable and uncertain renewable generation [1]-[2]. Poor system balancing performance, in a market environment, will usually lead to higher fluctuations in electricity prices, especially in the real-time (RT) markets of typical two-settlement markets [3]. In some markets, the frequency of price spikes is indeed increasing due to the changing generation portfolio [4]. In the US, RT markets are usually cleared every 5 minutes by the RT economic dispatch (RTED) tool, which considers a single time interval or multiple intervals [5]. When there are significant variations in load, interchange transactions or variable generation, the system could be in lack of flexible capacity to move from one RT dispatch point to the next. This challenge has motivated the introduction of flexible ramp capability (FRC) products, also termed flexiramp [6].

Current market practices in wholesale energy markets to compensate for such ramp shortages can be categorized as: increasing reserve margins, adding offset value to the fore-casted load, starting fast-start units, and keeping some additional units online [7]. In the US, some ISOs have launched markets for FRC. The practices mainly include the flexible ramping products implemented in California ISO and Midcontinent ISO [8]-[10]. Such products are co-optimized and cleared with energy and other ancillary services. Once obtained, the FRC capacities are withheld in the current period to meet the ramping requirements for future periods.

According to current practices, FRCs are mainly procured in RT markets. More specifically, fixed ramping requirements are embedded into the deterministic RTUC/RTED models [12]. Recent research shows that FRC can also be obtained through a stochastic problem which can achieve better performance [6]. However, RTED does not consider commitment of fast-start units, hence being less able to ensure the physical adequacy of FRCs. In [13], the procurement of FRC is considered in the RT unit commitment (RTUC). In this case, FRC requirements can trigger the commitment of fast-start units. Nevertheless, RTUC itself is not a market process. Procuring FRC by solving RTUC will consequently degrade into an out-of-market approach, which departs from the motivation of designing this product. Moreover, it is stated in [9] that because RTED has more accurate information than RTUC, the overall requirements of FRC should be less, and their procurement more efficient, compared to those in RTUC. It seems that how to properly clear this new product and meanwhile guarantee the adequate supply still remains as a controversial problem and unresolved.

In addition, only conventional units are currently entitled to provide FRCs to the system [14]. However, due to the increasing share of variable generation in the generation portfolio, especially wind power, the system will be facing increasing demands for FRC, while the resources that currently provide FRC are dispatched less frequently. If that is the case, holding to the current rules might lead to more FRC scarcity, which inspires us to consider the possibility of including variable generation producers, such as wind power producers (WPPs) as FRC providers. It is expected that, in this way, the usage of expensive fast-start units can be reduced. Actually, similar ideas have been presented by introducing WPPs into reserve and frequency regulation markets [15]-[17]. However, before implementing the mechanism, it has to be carefully examined whether it is economic to have low-marginal-cost technologies provide FRC products, and, how their capabilities of providing FRC should be defined considering their uncertain and volatile generation availability.

Following the recent trend of developing FRC markets, this paper explores the potential benefits of WPPs providing FRC services. As we will see in the following, WPPs bear the risk of being curtailed in the future periods while providing FRC services. In particular, WPPs have to reduce their output to provide upward FRC. Therefore, they only have comparative advantages providing FRCs when they can replace expensive fast-start units. To investigate the benefits and costs, a real-time operation process considering FRC products is modeled. The RTUC and RTED processes are formulated as a two-stage stochastic problem to consider the uncertain realization of wind and load. The first-stage problem represents the RTUC which does not include explicit FRC requirements but rather addresses the capability of obtaining adequate FRCs in the second-stage RTED process. The main contributions of this paper can be summarized as follows:

1) We propose to include WPPs as FRC providers and describe how they can be embedded into the FRC markets;

2) We propose a two-stage scheme for modeling the RT operation process including RTUC and RTED, in which the capabilities of WPPs to provide FRC are considered;

3) We carry out case studies to demonstrate the efficiency of including WPPs as FRC providers.

The remainder of this paper is organized as follows: in Section II, we describe the concepts of WPPs' providing FRC services; in Section III, we present the mathematical formulation of the two-stage RT operation and market process; cases studies are demonstrated in Section IV and V, followed by the conclusions in Section VI.

## II. WIND POWER PROVIDING FLEXIBLE RAMPING PRODUCT

FRC aims at reserving flexibility in a specific time slot for future use [18]. To satisfy the forecasted net load changes and cover the potential unexpected variations, the system needs to prepare adequate upward and downward FRCs see Fig. 1).

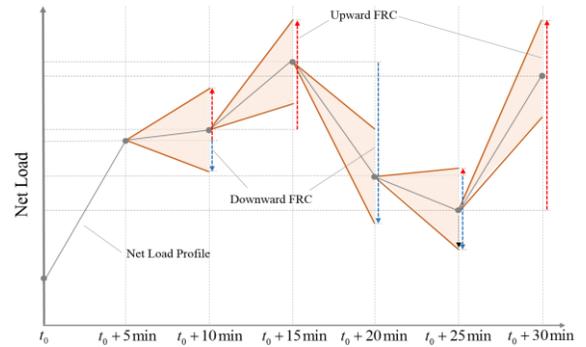

Fig. 1. Illustration of ramping capacity requirements in RT operation.

Generally, the FRCs are provided by online units, limited by their generation capacity, minimum output limit and ramp rate limit. The units that are willing to provide FRC services can submit bids including the above information to the system (currently only zero-price bids are allowed in CAISO and MISO). Then, market clearing results will be automatically decided by the RTED with FRC requirement constraints (single-period SCED is considered in this paper). When the overall ramp capacity is sufficient, these constraints will not bind, leading to zero FRC prices. When the ramp capacity becomes limited, the schedules of some units will need to be revised to release more FRCs. If the requirements still cannot be met, the system will be exposed to significant risks. Once the adjustable resources are unable to follow the change of load, the system will have to deploy regulation reserves by automatic generation control (AGC) on a continuous basis to maintain frequency. More severely, load shedding and wind spillage might happen.

Currently, WPPs are not considered as FRC providers. Instead, they have contributed to the increasing needs for flexibility. However, wind turbines are actually physically capable of providing FRC. The most widely used types of wind turbines, i.e., those based on double-fed induction generators and di-rect-drive permanent-magnet generators implemented with power electronic devices and advanced control technologies, are capable of adjusting their active power output rapidly, e.g., 0.05-0.25 p.u./s [19]. Comparatively, the availability of such advanced controls is even more advantageous in terms of flexibility than thermal generators.

Despite the facts above, WPPs are still excluded from FRC markets. This is largely due to two major concerns. Firstly, from the physical perspective, the availability of wind is uncertain, thus their capability of providing FRCs, though with adequate ramp rates, cannot be guaranteed because of the "real" capacity constraint that is dictated by the wind resource. Secondly, from the economics perspective, to provide upward FRC, the WPP has to reduce output to reserve a certain amount of upward capacity. Moreover, when FRCs are scheduled for a WPP, it is highly possible that wind spillage will happen during the following periods, as discussed in more detail later. In other words, the opportunity costs of WPPs, i.e. the lost energy payment, to provide FRCs is relatively high.

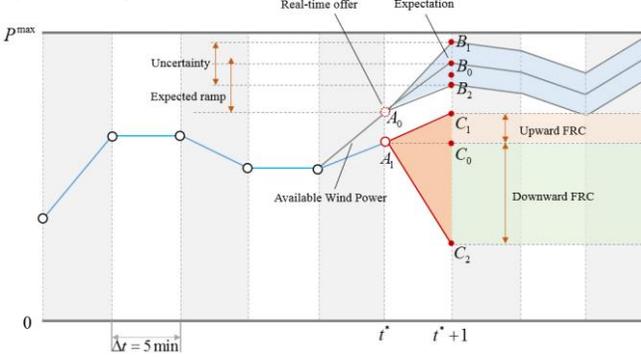

Fig. 2. Illustration of WPP's providing ramping capacity.

For the first concern, it should be noted that, though inherently uncertain, the forecast uncertainty is relatively low for the very short look-ahead horizon and WPPs are still capable of providing some amount of reliable FRCs. The key point of defining this capability is to differentiate between the passive variations from those that are controllable. Fig. 2 illustrates that when the RTED problem is solved for time $t^*$, a WPP submits its energy offer represented by point $A_0$, and its ramping capabilities. At that time, the expected output of the WPP at $(t^*+1)$ is $B_0$. But, due to uncertainty, the realized output, denoted by $B^*$ is highly possible to be located in $[B_1, B_2]$, which is a range that might be either wider or narrower than illustrated. Based on these assumptions, we can discuss the WPP's capability of providing FRC services below.

*1) Upward FRC*: if the WPP is scheduled to provide upward FRC, it has to deload by an amount of $(A_0-A_1)$. Then, for most of the situations, the WPP is able to provide $(A_0-A_1)$ of controllable upward ramping capability, except for some extremely low wind cases. To elaborate, comparing to letting the WPP passively ramp from $A_0$ to $B^*$, scheduling it at $A_1$ will give the system more controllable capacity, quantified as $(A_0-A_1)$, for the subsequent dispatch. Therefore, the capability of the WPP to provide upward FRC can be expressed by

$$r_{j,t,s}^{w,u} \triangleq \tilde{p}_{j,t,s}^{w} - p_{j,t,s}^{w} \tag{1}$$

$$r_{j,t,s}^{w,u} \leq ramp_j^w \cdot 5\min \tag{2}$$

*2) Downward FRC*: because we can always choose to curtail wind power in RTED if necessary, the capability of the WPP to provide downward FRC mainly depends on the wind power availability at $(t^*+1)$. To guarantee the deliverability of the downward FRC considering the forecasting uncertainty, we can limit the capacity by the expected wind power output subtracted by a safety margin, i.e., a quantile of probabilistic distribution at a specific confidence level (1-$\alpha$). Therefore, the capability of the WPP to provide downward FRC can be expressed by

$$r_{j,t,s}^{w,d} \leq \underline{p}_{j,t+1,s}^{w,\alpha} \tag{3}$$

$$r_{j,t,s}^{w,d} \leq ramp_j^w \cdot 5\min \tag{4}$$

For the concern on economic efficiency, WPPs usually have higher opportunity costs than conventional units to provide FRCs as WPPs normally get dispatched first due to their low marginal cost. Moreover, if upward FRCs provided by WPPs are not called for or downward FRCs are deployed, wind power spillage will happen. However, from the perspective of the entire system, the utilization of WPPs' ramp capability may help reduce the use of even more expensive thermal units. Hence, the WPPs may have comparative advantages conventional in some situations. For systems that are in lack of ramp resources, reserving FRCs preventively can be more economically efficient than taking corrective action afterwards. If the overall system costs can indeed be reduced by introducing WPPs in FRC markets, and they are properly compensated, including them as FRC providers is also economically viable.

### III. TWO-STAGE REAL-TIME UNIT COMMITMENT CONSIDERING FRC ADEQUACY

#### A. RT operation process overview

To cope with the uncertainties, multiple commitment and scheduling procedures with different lead times are conducted by system operators. Some of them are also market processes, while others are only for reliability concerns. Fig. 1 in [20] provides a good demonstration of the relationships among these procedures and the uncertainty levels they face. Fig. 3 in this paper is modified accordingly to show mainly the RT operation process, which we define to include the RTUC and RTED.

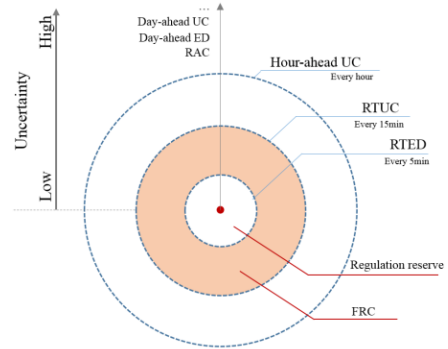

Fig. 3. RT operation process associated with the corresponding level of uncertainty (Modified based on [20]).

RTUC is also named as look-ahead unit commitment (LAC). It typically runs every 15 min with a study period of hours and resolution of 15 or 30 min. Hence, there are multiple RTED procedures between two RTUCs. RTUC only decides the commitment of fast-start units, as the commitment of slow-start units has been determined by previous scheduling processes. Therefore, the complexity of a RTUC is much less than the standard day-ahead unit commitment. The target of RTUC is to guarantee sufficient capacity and ramp capability. After the RTUC is solved, unit commitment status is held unchanged for a period of time. The system has to schedule sufficient FRCs from the online generation portfolio in subsequent RTEDs to keep itself capable of moving from one status to another, represented by the RTED results. Between two RTEDs, there is no scheduling process, and regulation reserves are deployed by AGC to continuously balance supply and demand while keeping the system safe from contingencies and large disturbances.

The RTUC problem is solved with the latest forecasts. It can be formulated either in a deterministic form with minimum FRC requirements, or, as a stochastic problem without explicit FRC constraints [6]. However, as mentioned above, in current markets the FRCs are usually procured in the RTED rather than RTUC. Namely, no matter what load profile is realized after the RTUC, the system will still purchase adequate FRCs at every step of the RTED. Therefore, it is more reasonable and practical to consider the RTED process embedded with FRC requirements in the RTUC problem formulation rather than only consider the FRC sufficiency at the RTUC level. In other words, following the real-world operation process, the RTUC is designed to guarantee the feasibility of subsequent RTEDs, and the RTEDs address the economic value of FRCs as well as further uncertainties by including FRC constraints. This leads to the formulation of a two-stage stochastic RTUC (TS-RTUC) problem. Such a model is presented in this paper to help reveal the effects of introducing WPPs as FRC providers.

### B. TS-RTUC problem formulation

The TS-RTUC problem is formulated to include a 15-min resolution RTUC problem in the upper level, and a 5-min resolution RTED problem in the lower level. The fast-start units can only be committed every 15 min. The objective function of the problem can be formulated as follows:

$$\min \sum_{t \in T} \sum_{i \in FG} u_{i,t} c_i^{SU} + \sum_{s \in S} \rho_s \left[ \frac{1}{12} \sum_{t \in T} \sum_{i \in G} C_i \left( p_{i,t,s}^g \right) + VOLL \cdot \sum_{t \in T} \Delta l_{t,s} + \omega \left( \Delta r_t^u + \Delta r_t^d \right) \right] \quad (5)$$

where, the startup and fuel costs are considered, as well as the penalty for unsatisfied loads and FRC shortages. The objective function is subject to:

*1) System balance constraints*

$$\sum_{i \in G} p_{i,t,s}^g + \sum_{j \in W} p_{j,t,s}^w = L_s + \Delta l_s \quad (6)$$

*2) Unit commitment constraints*

$$x_{i,t} - x_{i,t-1} = u_{i,t} - v_{i,t} \quad (7)$$

$$u_{i,t} = v_{i,t} = 0, \quad t \notin T^{15} \quad (8)$$

*3) Generation Capacity constraints of the units*

$$x_{i,t,s} P_i^{g,\max} - p_{i,t,s}^g - q_{i,t,s}^{g,u} - r_{i,t,s}^{g,u} \geq 0 \quad (9)$$

$$-x_{i,t,s} P_i^{g,\min} + p_{i,t,s}^g - q_{i,t,s}^{g,d} - r_{i,t,s}^{g,d} \geq 0 \quad (10)$$

*4) Regulation reserve requirements constraints*

$$\sum_{i \in G} q_{i,t,s}^{g,u} \geq Q^{u,req}, \quad \sum_{i \in G} q_{i,t,s}^{g,d} \geq Q^{d,req} \quad (11)$$

*5) Ramping constraints of the units*

$$p_{i,t+1,s}^g - p_{i,t,s}^g \leq x_{i,t,s} ramp_i \cdot 5\min + u_{i,t+1} P_i^{g,\min} \quad (12)$$

$$p_{i,t,s}^g - p_{i,t+1,s}^g \leq x_{i,t+1,s} ramp_i \cdot 5\min + v_{i,t} P_i^{g,\min} \quad (13)$$

$$0 \leq r_{i,t,s}^{g,u} \leq x_{i,t} ramp_i \cdot 5\min \quad (14)$$

$$0 \leq r_{i,t,s}^{g,d} \leq x_{i,t} ramp_i \cdot 5\min \quad (15)$$

$$r_{i,t,s}^{g,u} + q_{i,t,s}^{g,u} \leq x_{i,t} ramp_i \cdot 10\min \quad (16)$$

$$r_{i,t,s}^{g,d} + q_{i,t,s}^{g,d} \leq x_{i,t} ramp_i \cdot 10\min \quad (17)$$

$$r_{i,t,s}^{g,u} \leq \left(1 - x_{i,t-1,s}\right) SURR_i \cdot 5\min, \quad t \in T^{15} \quad (18)$$

$$r_{i,t-1,s}^d \leq x_{i,t,s} SDRR_i \cdot 5\min, \quad t \in T^{15} \quad (19)$$

*6) Non-negativity constraints*

$$p_{i,t,s}^g, p_{j,t,s}^w, \Delta l_s, q_{i,t,s}^{g,u}, q_{i,t,s}^{g,d} \geq 0 \quad (20)$$

Constraints (7)-(8) means that the unit commitments are decided at the upper level problem and fixed for all the scenarios at the lower level. Constraints (14)-(17) mean that the FRCs provided by the units are limited by their maximum ramp rates. Additionally, ramp capability within 10-min time scale is shared between the FRCs and scheduled regulation reserves according to [9]. Constraint (18) means that, when the unit is off at time (*t*-1) and time *t* is a RTUC point, it is capable to provide upward FRC up to its SURR. Similarly, constraint (19) means that when the unit is on previously, it can be shut down to provide downward FRC up to its SDRR. For some units, the SURR and SDRR are larger than their normal ramp rates.

However, some of the units, though categorized as fast-start units, cannot be completed started up or shut down within 5 min, meaning that after they are committed up or down, they might still be in their startup and shutdown processes at the subsequent RTED moment, being not controllable and can be regarded as contributing to the passive variation of load [22]. It is difficult to describe such processes precisely. But they can still be modeled based on realistic assumptions: (1) The fast-start units considered in the RTUC can be classified into two groups: type I and type II. Type I units can be started up or shut down within 5 min, including small-size gas turbines and hydro units, and do not need special treatment. Type II units include the other types of fast-start units which are assumed to take 15 min to start up or shut down, thus requiring special consideration; (2) The startup and shutdown trajectories are assumed to be linear, namely that they ramp up or down at a constant rate during the 15-min period. The assumptions about the startup and shutdown processes are a simplification of the actual situation, and can be furnished furtherly by adding more details. Based upon the assumptions, the startup process of the type II units can be represented by the following constraints, and the constraints of shutdown process are neglected for conciseness.

$$\overline{p}_{i,t,s}^{g,su} = u_{i,t} SURR_i \cdot 5\min + \left(1 - u_{i,t}\right) M, \quad t \in T^{15} \quad (21)$$

$$\underline{p}_{i,t,s}^{g,su} = u_{i,t} P_i^{g,\min} / 3 - \left(1 - u_{i,t}\right) M, \quad t \in T^{15} \quad (22)$$

$$\overline{p}_{i,t+1,s}^{g,su} = u_{i,t} SURR_i \cdot 10\min + \left(1 - u_{i,t}\right) M, \quad t \in T^{15} \quad (23)$$

$$\underline{p}_{i,t+1,s}^{g,su} = 2u_{i,t} P_i^{g,\min} / 3 - \left(1 - u_{i,t}\right) M, \quad t \in T^{15} \quad (24)$$

$$\overline{p}_{i,t+2,s}^{g,su} = u_{i,t} SURR_i \cdot 15\min + \left(1 - u_{i,t}\right) M, \quad t \in T^{15} \quad (25)$$

$$\underline{p}_{i,t+2,s}^{g,su} = u_{i,t} P_i^{g,\min} - \left(1 - u_{i,t}\right) M, \quad t \in T^{15} \quad (26)$$

The output of the units is limited by the auxiliary upper and lower bounds describing the power trajectory of startup and shutdown processes. Therefore, we have

$$\underline{p}_{i,t,s}^{g,su} \leq p_{i,t,s}^{g} \leq \overline{p}_{i,t,s}^{g,su} \quad (27)$$

$$\underline{p}_{i,t,s}^{g,sd} \leq p_{i,t,s}^{g} \leq \overline{p}_{i,t,s}^{g,sd} \quad (28)$$

Moreover, the passive ramping processes are nearly not controllable and should be considered in the load variation. The amount of passive ramping of a unit can be represented by

$$\Delta p_{i,t,s}^{g,su}, \Delta p_{i,t+1,s}^{g,su}, \Delta p_{i,t+2,s}^{g,su} = u_{i,t} P_i^{g,\min}/3, \ t \in T^{15}, i \in G^{II} \quad (29)$$

$$\Delta p_{i,t,s}^{g,sd}, \Delta p_{i,t+1,s}^{g,sd}, \Delta p_{i,t+2,s}^{g,sd} = -v_{i,t} P_i^{g,\min}/3, \ t \in T^{15}, i \in G^{II} \quad (30)$$

They are approximations of the passive ramping, the extra portion caused by ramping up to above $P^{\min}$ or ramping down from above $P^{\min}$ is ignored. In addition, type II units during their startup and shutdown processes are unable to provide flexibility, therefore

$$r_{i,t,s}^{g,u}, r_{i,t+1,s}^{g,u}, r_{i,t+2,s}^{g,u} \leq (1 - u_{i,t} - v_{i,t})M, \ t \in T^{15}, g \in G^{II} \quad (31)$$

$$r_{i,t,s}^{g,d}, r_{i,t+1,s}^{g,d}, r_{i,t+1,s}^{g,d} \leq (1 - u_{i,t} - v_{i,t})M, \ t \in T^{15}, g \in G^{II} \quad (32)$$

At last, for every possible scenario, we have to guarantee that adequate FRCs can be procured in the RTED process. That is to say, we are not trying to cover all the possibilities in the RTUC problem. Instead, sufficient flexibility is guaranteed for each of the limited scenarios to cover further uncertainties. It is somehow closer to the real-world process as discussed above. The FRC constraints can be formulated as follows:

$$\sum_{i \in G} r_{i,t,s}^{g,u} \geq \left[ \hat{L}_{t+1,s} - \sum_{i \in G^{II}} \left( \Delta p_{i,t+1,s}^{g,su} - \Delta p_{i,t+1,s}^{g,sd} \right) - \sum_{j \in W} \hat{p}_{j,t+1,s}^{w} \right] \\ - \left( L_{t,s} - \Delta l_{t,s} - \sum_{j \in W} \tilde{p}_{j,t,s}^{w} \right) + \xi^{u} \quad (33)$$

$$\sum_{i \in G} r_{i,t,s}^{g,d} \geq -\left[ \hat{L}_{t+1,s} - \sum_{i \in G^{II}} \left( \Delta p_{i,t+1,s}^{g,su} - \Delta p_{i,t+1,s}^{g,sd} \right) - \sum_{j \in W} \hat{p}_{j,t+1,s}^{w} \right] \\ + \left( L_{t,s} - \Delta l_{t,s} - \sum_{j \in W} \tilde{p}_{j,t,s}^{w} \right) + \xi^{d} \quad (34)$$

The first two terms in (33)-(34) represent the expected passive variation of the loads, consisting of changes in loads, wind power and uncontrollable ramping of the type II units. The last term represents an extra margin, decided mainly based upon the distribution and specified confidence level of the load and wind power forecast error [21]. Combining (33) and (34), it is required that the system is able to cover a specific range of load variation by obtaining FRCs.

The objective function (5) subject to constraints (6)-(34) and those describing the power trajectories of startup and shutdown processes forms the optimization problem for the TS-RTUC.

### C. Consider WPPs as FRC providers

The above model can be easily modified to include WPPs as FRC providers. To elaborate, the constraints (1)-(4) are included, and the constraints (33)-(34) are modified as

$$\sum_{i \in G} r_{i,t,s}^{g,u} + \sum_{j \in W} r_{i,t,s}^{w,u} \geq \left[ \hat{L}_{t+1,s} - \sum_{i \in G^{II}} \left( \Delta p_{i,t+1,s}^{g,su} - \Delta p_{i,t+1,s}^{g,sd} \right) \\ - \sum_{j \in W} \hat{p}_{j,t+1,s}^{w} \right] - \left( L_{t,s} - \Delta l_{t,s} - \sum_{j \in W} \tilde{p}_{j,t,s}^{w} \right) + \xi^{u} \quad (35)$$

$$\sum_{i \in G} r_{i,t,s}^{g,d} + \sum_{j \in W} r_{i,t,s}^{w,d} \geq -\left[ \hat{L}_{t+1,s} - \sum_{i \in G^{II}} \left( \Delta p_{i,t+1,s}^{g,su} - \Delta p_{i,t+1,s}^{g,sd} \right) \\ - \sum_{j \in W} \hat{p}_{j,t+1,s}^{w} \right] + \left( L_{t,s} - \Delta l_{t,s} - \sum_{j \in W} \tilde{p}_{j,t,s}^{w} \right) + \xi^{d} \quad (36)$$

Furthermore, some modifications are needed to consider the economic effects of having WPPs to provide FRC services. Currently, the FRC providers are compensated by FRC credits, which depend on the multipliers of constraints (35)-(36) derived by a RTED problem. Because FRCs are co-optimized with energy products, the opportunity cost of providing FRCs and not producing energy (or producing at marginal costs higher than the clearing energy prices) are automatically compensated by the FRC prices [3]. However, such credits can only cover the potential costs for the current period. When some of the FRCs are procured from the WPPs, the system might have even less flexibility provided by the traditionally defined. Consequently, the possibility of wind power spillage in future periods might increase. This risk is not reflected by the compensation as this is based solely upon the opportunity costs of the current period. Therefore, one can argue that extra compensation should be given to the WPPs for providing FRC services and bearing underlying risks. From the system's perspective, this extra compensation stands for the potential subsequent avoidance of wind power spillage caused by scheduling WPPS as FRC sources. Based on the above discussion, the objective function can be modified as

$$\min \sum_{t \in T} \sum_{i \in FG} u_{i,t} c_i^{SU} + \sum_{s \in S} \rho_s \left[ \frac{1}{12} \sum_{t \in T} \sum_{i \in G} C_i \left( p_{i,t,s}^{g} \right) \\ + VOLL \cdot \sum_{t \in T} \Delta l_{t,s} + \omega \left( \Delta r_t^{u} + \Delta r_t^{d} \right) + \sum_{t \in T} \sum_{j \in W} \left( \pi^{u} r_{j,t,s}^{w,u} + \pi^{d} r_{j,t,s}^{w,d} \right) \right] \quad (37)$$

A reasonable evaluation on the values of $\pi^u$ and $\pi^d$ is not straightforward and out of scope of this paper. It can be obtained by statistical analysis. Or, a more market-based approach would be to allow the WPPs to bid at a non-zero price and let themselves decide the proper level of compensation. But it should be noted that including such compensation is relatively neutral to the model. If the integrated costs of using WPP-provided FRC services are too high, the optimization problem will just automatically choose not to use them. In the case study section, we will explore the effects of different compensation factors.

### IV. CASE STUDY: MODEL AND ASSUMPTIONS

Following [13], we use a system including five generators in our studies. The generator data (Table I) is modified based on [13]. Among the generators, G1 is a base-load unit, and G2-G5 are fast-start units with increasing fuel costs, which are assumed to be linear with respect to power output, i.e., $C_i(p) = a_i p + b_i$. According to the definition in Section III, G2-G4 are type II units, while G5 is a type I. A wind farm with an installed capacity of 150 MW is also included in the system. The maximum ramping rate of the wind farm is set as 20 MW/min.

In practice, the RTUC problem is solved in a rolling manner. However, we only consider one hour and assume that no new forecast information is provided within the hour. Hence the

RTUC problem is solved without rolling and considering the one hour study period, which consists of 4 RTUC points and 12 RTED points.

The extra upward and downward FRC margins, i.e., $\xi^+$ and $\xi^-$ are assumed to be linear functions of system-level forecasts of load and wind power, which can be expressed as follows:

$$\xi^u = \eta^{l,u} \hat{L}_{t+1} + \eta^{w,u} \sum_{j \in W} \hat{p}^w_{j,t+1,s} \quad (38)$$

$$\xi^d = \eta^{l,d} \hat{L}_{t+1} + \eta^{w,d} \sum_{j \in W} \hat{p}^w_{j,t+1,s} \quad (39)$$

In the cases, if not specially noted, the parameters are set as $\eta^{d,u} = \eta^{l,d} = 0.01$, and $\eta^{w,u} = \eta^{w,d} = 0.05$. For real systems, the determination of a proper margin can be tuned based on statistical analysis of practical data. The load profile used in this study is obtained from BPA [23] and scaled down to fit the test system. The wind power data is modified based on [24]. For each TS-RTUC, 5000 scenarios are generated respectively for load and wind power, and are both reduced to 10 scenarios using the scenario reduction technique [25], making up 100 scenarios in total. While generating scenarios, load is assumed to follow a normal distribution and wind power is assumed to follow a beta distribution. The standard deviation of load and wind power forecast are assumed to be 0.5% and 1% of the expectation. In addition, it is assumed that the WPP always offers its expected wind power to the system, i.e. $\tilde{P}^w_{j,t,s} = \hat{P}^w_{j,t,s}$. For unsatisfied load and FRC requirements, the *VOLL* is set to be $4000/MWh, and the FRC scarcity price is set to be $1000/MWh. The optimization problems are implemented in a Matlab environment and solved by Gurobi Optimizer.

TABLE I
GENERATOR CHARACTERISTICS

| Gen | $P^{max}$ (MW) | $P^{min}$ (MW) | ramp (MW/min) | SURR (MW/min) | SDRR (MW/min) | a ($/MWh) | b ($/h) | $c^{SU}$ ($) |
|---|---|---|---|---|---|---|---|---|
| G1 | 300 | 300 | 0 | 0 | 0 | 10 | 0 | 0 |
| G2 | 150 | 50 | 3 | 4 | 4 | 20 | 300 | 300 |
| G3 | 200 | 50 | 3 | 4 | 4 | 40 | 300 | 600 |
| G4 | 150 | 50 | 3 | 5 | 5 | 60 | 300 | 900 |
| G5 | 100 | 10 | 6 | 6 | 6 | 120 | 0 | 0 |

## V. NUMERICAL RESULTS

### A. Single-hour cases

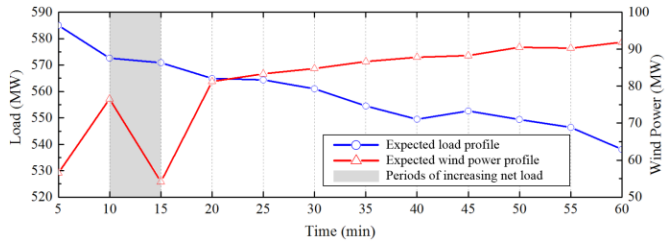

Fig. 4. Expected load and wind power profile of a specific study case.

First, several single-hour cases are carried out to demonstrate how a WPP can act as a FRC provider and its effects. In Fig. 4, the expected load and wind power profiles during a certain hour are presented. Based on the profiles, the TS-RTUC problem with and without the WPP's providing FRC are solved respectively for comparison. For conciseness, the former case is denoted by "W-FRC", and the latter case "NW-FRC".

The unit commitment results as shown in Table II, which indicate that by introducing the WPP as a FRC provider, the commitment of the very expensive unit G5 is avoided. And, even though extra compensation is given to the WPP for providing FRC ($\pi^+ = \pi^- = \$20$/MWh), the expected overall cost is still reduced by 5.5% compared to the NW-FRC case. Fig. 5 presents the average portfolio of FRC provision under the two cases. The comparison clearly presents the way that the WPP contributes to the system flexibility. To elaborate, for the first half hour, the FRCs provided by the WPP make it unnecessary to start G5 in the W-FRC case. For the second half hour, especially at 45 and 50 min, detailed results show that the system experiences shortage of FRCs under some rare scenarios. Because of the very low probability, G5 is not started for economics consideration in the NW-FRC case. But, in the W-FRC case, the ramping capability requirements can still be satisfied in all scenarios.

TABLE II
UNIT COMMITMENT RESULT COMPARISON

| Gen | Unit Commitment (W-FRC/NW-FRC) Time (min) | | | | | | | | | | | |
|---|---|---|---|---|---|---|---|---|---|---|---|---|
| | 5 | 10 | 15 | 20 | 25 | 30 | 35 | 40 | 45 | 50 | 55 | 60 |
| G1 | 1/1 | 1/1 | 1/1 | 1/1 | 1/1 | 1/1 | 1/1 | 1/1 | 1/1 | 1/1 | 1/1 | 1/1 |
| G2 | 1/1 | 1/1 | 1/1 | 1/1 | 1/1 | 1/1 | 1/1 | 1/1 | 1/1 | 1/1 | 1/1 | 1/1 |
| G3 | 1/1 | 1/1 | 1/1 | 1/1 | 1/1 | 1/1 | 1/1 | 1/1 | 1/1 | 1/1 | 1/1 | 1/1 |
| G4 | 0/0 | 0/0 | 0/0 | 0/0 | 0/0 | 0/0 | 0/0 | 0/0 | 0/0 | 0/0 | 0/0 | 0/0 |
| G5 | 0/1 | 0/1 | 0/1 | 0/1 | 0/1 | 0/1 | 0/0 | 0/0 | 0/0 | 0/0 | 0/0 | 0/0 |
| Expected generation costs (W-FRC/NW-FRC): $9740.01 / $11952.85 | | | | | | | | | | | | |

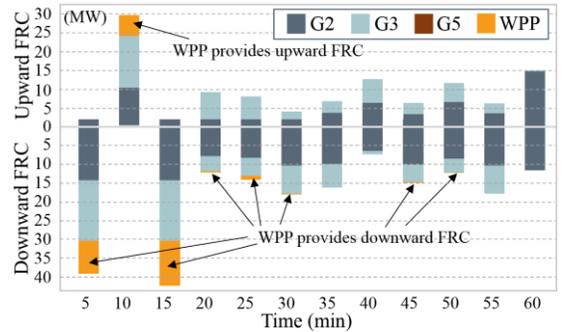

(a) W-FRC case

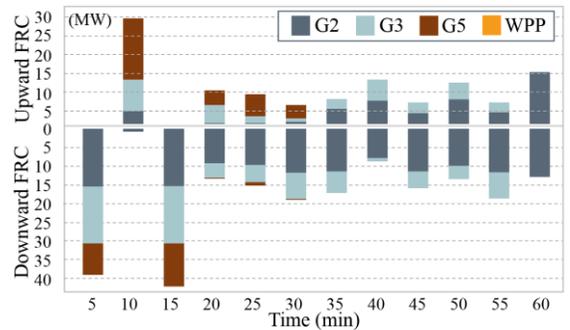

(b) NW-FRC cases

Fig. 5. Expected FRC provision portfolio with and without WPP-provided FRCs.

The expected wind power curtailments are 0.6072 MWh for the W-FRC case and 0.1914 MWh for the NW-FRC case. The increase in wind power curtailment is because the WPP reduces its output to provide upward FRCs at 5 min in the W-FRC case. Despite of the loss of wind power, the WPP is expected to earn a $56.24 credit for the service it provides. As long as the system is better-off and the WPP can be reasonably compensated, a moderate amount of wind power curtailment is acceptable. Table III provides a brief summary of three other typical cases. The results indicate that when the system is not in lack of ramping capability, the WPP will not be chosen to provide FRCs (case 1). This is because its relatively higher opportunity costs and extra compensation. As the ramping requirement increases (case 2 and case 3), the value of having WPP to provide FRCs also increases. This trend implies that, though not necessary now, in the future, when the system is deeply penetrated by wind power and other variable generation technologies, the idea of including WPPs as FRC providers might be a better choice than only using the conventional units.

TABLE III
SUMMARY OF UNIT COMMIT RESULTS FOR THREE TYPICAL CASES

| Case | Total Cost ($) | | Cost reduction (%) | G5 running hour (h) | | Wind power curtailed (MWh) | | WPP FRC credits ($) |
|---|---|---|---|---|---|---|---|---|
| | W-FRS | NW-FRS | | W-FRS | NW-FRS | W-FRS | NW-FRS | |
| 1 | 16307.5 | 16307.5 | 0.0 | 0.00 | 0.00 | 0.00 | 0.00 | 0.00 |
| 2 | 17860.4 | 18437.5 | 3.1 | 0.00 | 0.25 | 0.10 | 0.00 | 7.95 |
| 3 | 10876.5 | 12169.7 | 10.6 | 0.00 | 0.75 | 0.42 | 0.01 | 29.45 |

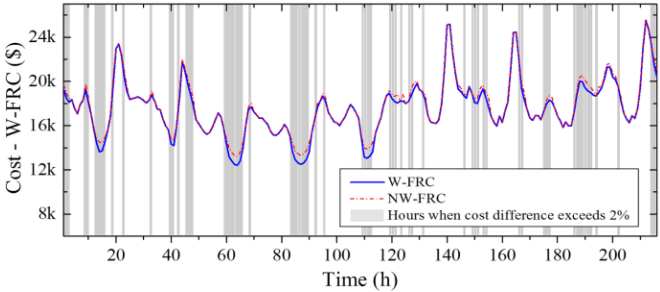
Fig. 6. Expected generation costs profile.

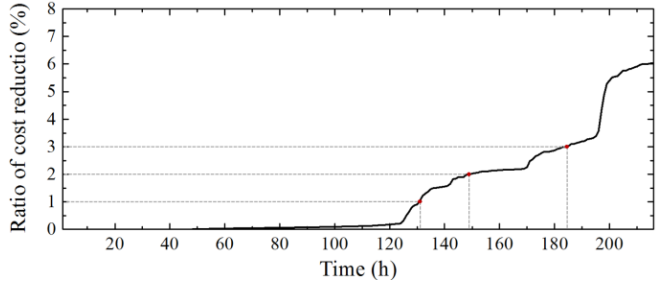
Fig. 7. Ratio of expected generation costs reduction sorted by the magnitude.

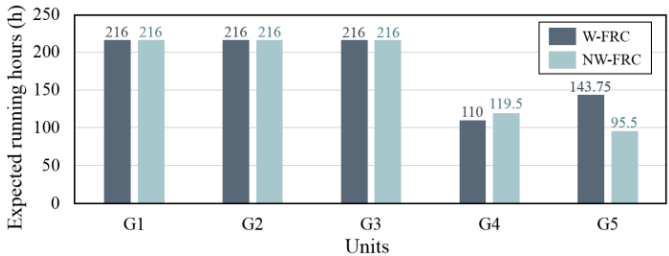
Fig. 8. Expected overall running hours of the units under the two cases.

### B. Multi-hour cases

In this subsection, cases including 216 consecutive hours, i.e., 9 days, are studied to further illustrate the effects of WPP-provided FRCs. When the parameters are the same as in the cases in the previous subsection, the expected generation costs for the W-FRC and NW-FRC cases are shown in Fig. 6. Fig. 7 presents the ratio of cost reduction sorted by the magnitude. The results indicate that costs are reduced by more than 1% for over 40% of the period, and 2% for over 30% of the period, and that the cost is never higher in the W-FRC case. Fig. 8 shows the expected overall running hours of the units under the two cases. As expected, the running hours of the expensive unit G5 are remarkably reduced in the W-FRC case. Instead, the running hours of relatively slow unit G4 are slightly increased. In Fig. 9, the expected wind power curtailment under both cases and the FRC credits earned under the W-FRC case are demonstrated. It is illustrated that, for most hours, more wind power will be curtailed in W-FRC case because the WPP is scheduled to provide upward FRC service. However, the amount is moderate, and it is compensated by the FRC credits. Additionally, for some hours the WPP does not have to be curtailed to earn credits because they can just provide downward FRC service.

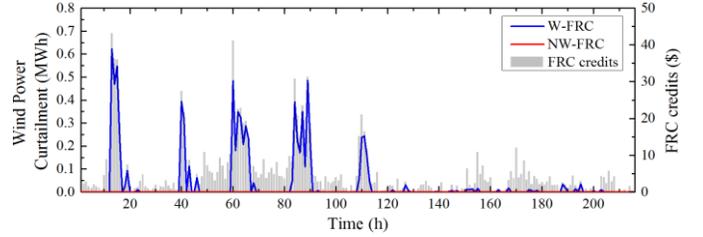
Fig. 9. Expected wind power curtailment and FRC credits.

TABLE IV
SENSITIVITY ANALYSIS ON THE COMPENSATION FACTORS

| $\pi^u, \pi^d$ ($/MWh) | Expected Ratio of cost reduction* (%) | Expected increase* of wind power spillage (MWh) | Expected WPP FRC credits ($) | Expected reduction* of running hours of G5 (h) |
|---|---|---|---|---|
| 0 | 1.17 | 15.19 | 566.57 | 56.25 |
| 10 | 1.15 | 12.05 | 949.85 | 54.25 |
| 20 | 1.14 | 8.43 | 1071.41 | 54.25 |
| 40 | 1.13 | 5.65 | 1313.16 | 53.75 |
| 80 | 1.10 | 5.64 | 1510.33 | 56.75 |
| 100 | 1.09 | 5.64 | 1282.51 | 55.75 |

*Increase and reduction are calculated as by comparing W-FRC case to NW-FRC case

TABLE V
SENSITIVITY ANALYSIS ON THE LEVEL OF UNCERTAINTY IN RTED PROCESS

| $\eta^{l,u}, \eta^{l,d}$ | $\eta^{w,u}, \eta^{w,d}$ | Expected Ratio of cost reduction* (%) | Expected increase* of wind power spillage (MWh) | Expected WPP FRC credits ($) | Expected reduction* of running hours of G5 (h) |
|---|---|---|---|---|---|
| 0.012 | 0.060 | 2.59 | 13.41 | 1852.42 | 85.75 |
| 0.010 | 0.050 | 1.14 | 8.43 | 1071.41 | 54.25 |
| 0.008 | 0.040 | 1.07 | 5.55 | 728.14 | 45.75 |
| 0.006 | 0.030 | 1.03 | 4.42 | 470.14 | 12.50 |
| 0.004 | 0.020 | 0.97 | 1.73 | 222.84 | 8.75 |
| 0.002 | 0.010 | 0.89 | 1.41 | 149.68 | 6.25 |

*Increase and reduction are calculated as by comparing W-FRC case to NW-FRC case

At last, some sensitivity analyses are carried out. First, the effects of different compensation factors for WPP-provided FRCs are examined. Table IV provides a summary of the test results, which indicate that, by providing an appropriate amount of compensation, the system cost is reduced and the WPP may also be better-off. Note that the WPP receives a significant FRC credit even when the compensation factor is zero. Moreover, even when the compensation level is very high, i.e. more than $40/MWh, there is still a significant reduction in the overall generation cost and running hours of G5, meaning that even when the FRC services provided by WPPs are highly priced, it is still more economic to use them to cover the FRC

requirements than to keep the expensive fast-start units stand by. Moreover, as mentioned in Section III, instead of fixing the compensation factor at system-level, we can try to allow the WPPs to submit non-zero bids for providing FRC services. It should be noted that the total revenue earned from providing energy by WPPs are not included in this study because we do not model the entire market settlement process where the WPPs may participate, especially the day-ahead energy markets. However, considering the two-settlement mechanism, fluctuating RT prices are not favorable. Therefore we assume FRC adequacy and reduced RT price volatility also a merit for WPPs. Further market simulations will be included in future studies.

Table V presents the comparison of results under different levels of uncertainties considered in the RTED process. The trends imply that, the value of WPPs' providing FRCs is more significant in situations with larger uncertainties. In lower uncertainty cases the WPP is not frequently scheduled to provide FRCs because of the underlying opportunity costs and extra compensation.

## VI. Conclusions

In this paper, the mechanism and possibility of including WPPs as FRC providers have been explored. Though inherently uncertain and variable, WPPs are capable of providing FRCs with good performance with the help of the very flexible power electronics devices embedded in modern wind turbines. However, it is also noticeable that the WPPs are not advantageous in providing such capacity-based services because of their low marginal costs and the increased likelihood of being curtailed. However, as the needs for FRCs significantly increase, the participation of WPPs in this market is expected to significantly reduce the commitment of expensive fast-start units, and reduce the risk of insufficient ramping capacity in the RT operation processes.

To reveal the effects, a TS-RTUC problem considering FRC adequacy is formulated. In this model, the upper level represents the RTUC of which the target is to guarantee that adequate FRCs can be obtained in the subsequent RTED processes, considering load and wind power uncertainties described by scenarios. The lower level is the RTED problem which co-optimizes energy, reserve and FRC products. Case studies based on this model demonstrate that WPPs' providing FRC services give significant reductions in overall system costs and WPPs may also earn significant FRC credits depending on the compensation scheme. The model can be extended and applied to markets where look-ahead (multi-period) ED tool is used.

In the future, we plan to extend this study to real systems with current or planned future wind power penetration levels to test the practical feasibility and economic efficiency of using wind power as FRC sources.